\documentclass[11pt]{article}
\usepackage{a4,amsthm,amssymb,amsmath,graphicx}

\theoremstyle{plain}
\newtheorem{thm}{Theorem}
\newtheorem{pro}{Proposition}

\theoremstyle{remark}
\newtheorem*{rem}{Remark}

\newcounter{abbildung}
\def\C{{\mathbb{C}}}
\def\epsilon{\varepsilon}
\def\e{\varepsilon}
\def\phi{\varphi}
\def\F{{\cal{F}}}
\def\B{{\cal{B}}}
\def\P{{\cal{P}}}
\def\L{{\cal{L}}}
\def\G{{\cal{G}}}
\def\hatP{\widehat{\cal{P}}}
\def\hatL{\widehat{\cal{L}}}
\def\tilP{\widetilde{\cal{P}}}
\def\tilL{\widetilde{\cal{L}}}
\def\tilV{\widetilde{V}}
\def\barP{\overline{\cal{P}}}
\def\barL{\overline{\cal{L}}}
\def\Qua{{\cal{Q}}}
\def\Se{{\cal{S}}}

\def\leer{\varnothing}
\def\bs{\backslash}
\def\Aut{\mathrm{Aut}}
\def\PGL{\mathrm{PGL}}
\def\iprod{\;\rule{0.04em}{0.6em}\rule{0.56em}{0.04em}\;}
\def\Wedge{\mbox{\Large$\wedge$}}
\def\chara{{\mathop{\rm char}\nolimits}\,}

\begin{document}
\title{On automorphisms of flag spaces%
\thanks{Scientific and technological cooperation Italy -- Austria
1998/2000, project no.~10.}}
\author{Hans Havlicek \and Klaus List \and Corrado Zanella}
\date{}
\maketitle

\begin{abstract}
We show that the automorphisms of the flag space associated with a
$3$-dimensional projective space can be characterized as bijections
preserving a certain binary relation on the set of flags in both
directions. From this we derive that there are no other automorphisms
of the flag space than those coming from collineations and dualities
of the underlying projective space. Further, for a commutative ground
field, we discuss the corresponding flag variety and characterize its
group of automorphic collineations.


\end{abstract}

\section{Introduction}\label{sec:intro}

The aim of the present article is to determine all automorphisms of
the \emph{flag space} associated with a $3$-dimensional projective
space $(\P, \L)$; cf.\ Section \ref{sec:flag}. Such an automorphism
is a bijection on the set $\F$ of flags that preserves pencils of
flags in both directions. However, we adopt a slightly different
point of view: Two flags are called \emph{related} ($\sim$) if they
differ in at most one of their three components. Now we ask for all
bijections $\alpha:\F\to \F$ such that
\begin{equation}\label{eq:0}
  \Phi\sim\Psi \Leftrightarrow \Phi^\alpha\sim\Psi^\alpha
  \mbox{ for all }\Phi,\,\Psi\in\F.
\end{equation}
Clearly, each collineation and each duality of $(\P, \L)$, via its
action on the set of flags, is a solution of (\ref{eq:0}). It will be
established in Section \ref{sec:pluecker}, that there are no other
solutions. Since the pencils of flags are exactly the maximal sets of
mutually related flags, this solves at the same time the problem to
find all automorphisms of the flag space.

Our result may also be seen as a characterization of the group of
collineations and dualities of a $3$-dimensional projective space
under a \emph{mild hypothesis} \cite{benz-92}, \cite{benz-94}. See
also \cite{pamb-00} for the logical background of such
characterizations.

In Section \ref{sec:variety} we turn to the classical point model of
$\F$, i.e. a \emph{flag variety}. (It is necessary here to assume
that the ground field is commutative.) We sketch a coordinate-free
approach using tools from multilinear algebra.  So, finally, we get
an intrinsic characterization of the group of collineations fixing
such a flag variety.

The flag variety associated with the $n$-dimensional projective space
over the complex numbers $\C$ has been discussed by \textsc{W.~Burau}
in \cite{bura-54}, \cite{bura-58}, \cite{bura-67}, and
\cite{bura-77}. One of the crucial tools in those papers is that this
flag variety yields an irreducible representation of the projective
group $\PGL(n+1,\C)$. If $\C$ is replaced with an arbitrary
commutative ground field $K$ then one still gets a representation of
$\PGL(n+1,K)$ as a group of projective collineations fixing the
associated flag variety. However, this representation is not
necessarily irreducible. So, in general the connection to irreducible
representations of linear groups is lost. For example, if the ground
field $K$ has characteristic $3$, then there is an invariant point in
the ambient space of the variety representing the flags of a
projective plane over $K$ \cite{thas+m-99}. Let us just mention here
that also for $n=3$ the representation of $\PGL(4,K)$ turns out to be
reducible in case of characteristic $3$, but this will be discussed
elsewhere.

\section{The flag space}\label{sec:flag}
Let $(\P, \L)$ be a 3-dimensional projective space with point set
$\P$ and line set $\L$. The subspaces of $(\P, \L)$ are considered
as sets of points. We shall not distinguish between a point
$Q\in \P$ and the subspace $\{Q\} \subset \P$. The sign $\vee$
is used to denote the join of projective subspaces.

Recall that a \emph{flag} is a triple $(P, g
, \e )$ consisting of a point $P$, a line $g$, and a plane $\e$ such
that $P\in g \subset \e$. We put $\F$ for the set of all flags of
$(\P, \L)$. Two flags $\Phi, \Psi \in \F$ are called \emph{related}
($\Phi \sim \Psi$) if they differ in at most one of their components.
We say that $\Phi$ and $\Psi$ are \emph{adjacent} if they are related
and distinct.

It is easy to show that $(\F, \sim)$ is a \emph{Pl\"ucker space} in
the sense of \textsc{W.~Benz} \cite[p.~199]{benz-92}: Clearly, the
relation $\sim$ is reflexive and symmetric. In order to show that
$(\F, \sim )$ is connected we consider two arbitrary flags $(P,g,\e)$
and $(P',g',\e')$. Then there is a line $h$ skew to $g$ and $g'$, a
point $Q\in h$ that is neither in $\e$ nor in $\e'$ and a plane
$\phi\supset h$ that contains neither $P$ nor $P'$. We infer that the
four points $P_0:=P$, $P_1:=g\cap\phi$, $P_2:=h\cap\e$, and $P_3:=Q$
form a tetrahedron (figure~\ref{abb1}). Put $g_{ij}$ for the edge
joining $P_i$ with $P_j$ and $\e_i$ for the face opposite to $P_i$.
Then $(P,g,\e)= (P_0,g_{01},\e_3)\sim (P_1,g_{01},\e_3)\sim
(P_1,g_{12},\e_3)\sim (P_1,g_{12},\e_0)\sim (P_2,g_{12},\e_0)\sim
(P_2,g_{23},\e_0)\sim (P_3,g_{23},\e_0)= (Q,h,\phi)$. Similarly,
$(Q,h,\phi)$ and $(P',g',\e')$ give rise to a tetrahedron, whence the
assertion holds.

{\unitlength1cm
    \begin{center}
      \begin{picture}(7,3.5)
         \put(2,0)
         {\includegraphics[width=3.5cm]{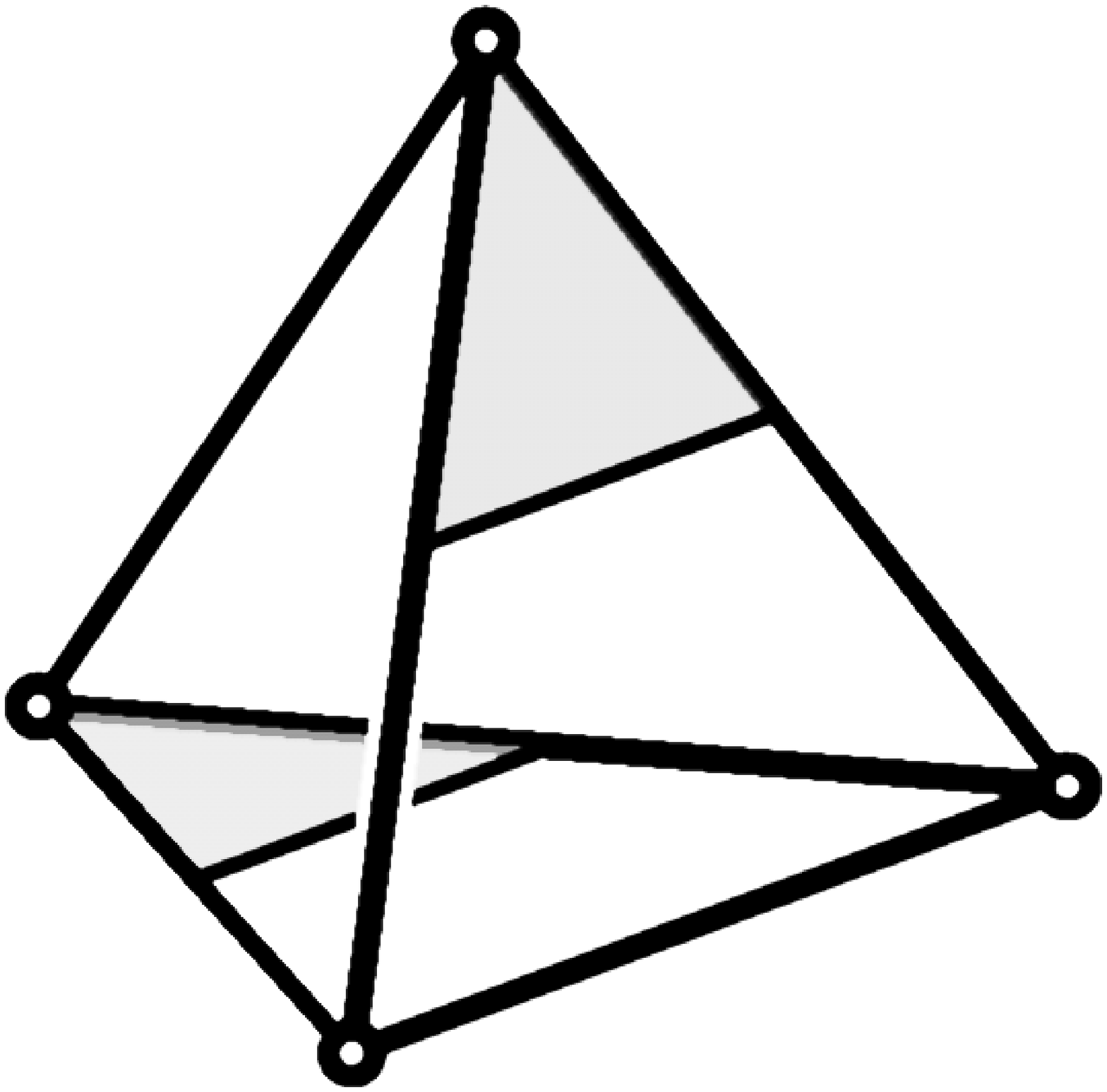}
         }
         \put(0.85 ,0.8){$P_0=P$}
         \put(5.45 ,0.6){$P_2$}
         \put(3.3 ,-0.25){$P_1$}
         \put(3.7 ,3.3){$P_3=Q$}
         \put(2.4 ,0.3){$g$}
         \put(4.6 ,2.1){$h$}
         \put(2.6 ,0.9){$\e$}
         \put(3.65,2.35){$\phi$}
       \end{picture}
       \vspace{.3cm}
       {\refstepcounter{abbildung}\label{abb1}
         \centerline{Fig.\ \ref{abb1}.}}
    \end{center}
 }

Let $P \in \P$ be a point. Then $\F [P]\subset \F$ stands for all
flags with point component $P$ and arbitrary other components. Given
a line $g\in \L$ or a plane $\e \subset \P$ then $\F [g]$ and $\F[\e
]$ are defined similarly. In addition we put $\F [P,g] := \F [P] \cap
\F [g]$, $\F [P,\e ] := \F [P] \cap \F [\e ]$, and $\F [g,\e]  := \F
[g] \cap \F [\e ]$.

The set $\F$ has three families of distinguished subsets, namely the
set $\B_i$ of \emph{pencils of type} $i$, where $i\in \{0,1,2\}$: A
pencil of type $0,1,2$ is a non-empty set of the form
\begin{equation}
  \label{eq:1}
  \F [g,\e],\; \F [P,\e ],\; F [P,g],
\end{equation}
respectively, where $P\in \P$, $g\in \L$, and $\e \subset \P$ is a
plane (figure~\ref{abb2}, \ref{abb3}, \ref{abb4}). Observe that
$\F [g , \e] \ne \leer$ is equivalent to $g\subset \e$ etc.\

For all flags of a pencil of type $i$ the components of (projective)
dimension $\ne i$ are the same, whence a pencil of type $i$ and a
pencil of type $j\ne i$ cannot coincide. We put $\B:= \B_0 \cup \B_1
\cup \B_2$ for the set of all \emph{pencils}. Each pencil contains as
many flags as there are points on a line.
{\unitlength1cm
      \begin{center}
      \begin{minipage}[t]{4.0cm}
         \begin{picture}(4.0,4.0)
         \put(0,0.7)
         {\includegraphics[width=4.0cm]{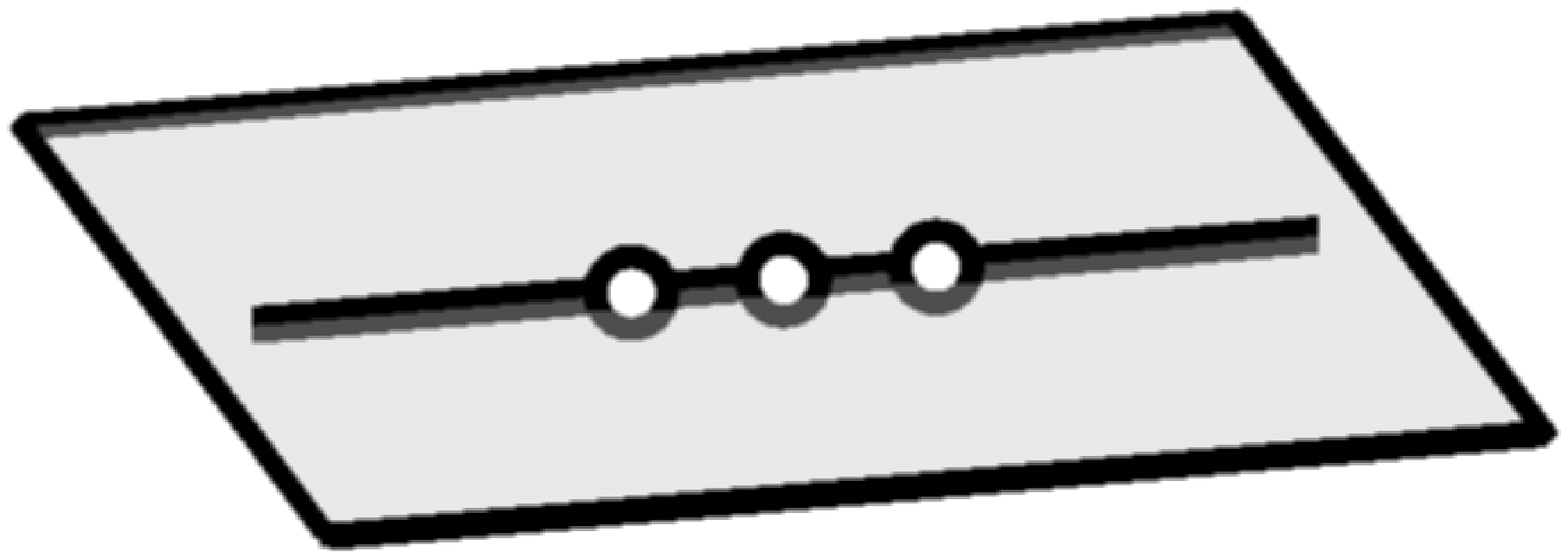}}
         \put(0.9 ,0.95){$\e$}
         \put(2.85 ,1.7){$g$}
         \end{picture}
         {\refstepcounter{abbildung}\label{abb2}
          \centerline{Fig.\ \ref{abb2}.}}
         \end{minipage}
      \begin{minipage}[t]{4.0cm}
         \begin{picture}(4.0,4.0)
         \put(0,0.7)
         {\includegraphics[width=4.0cm]{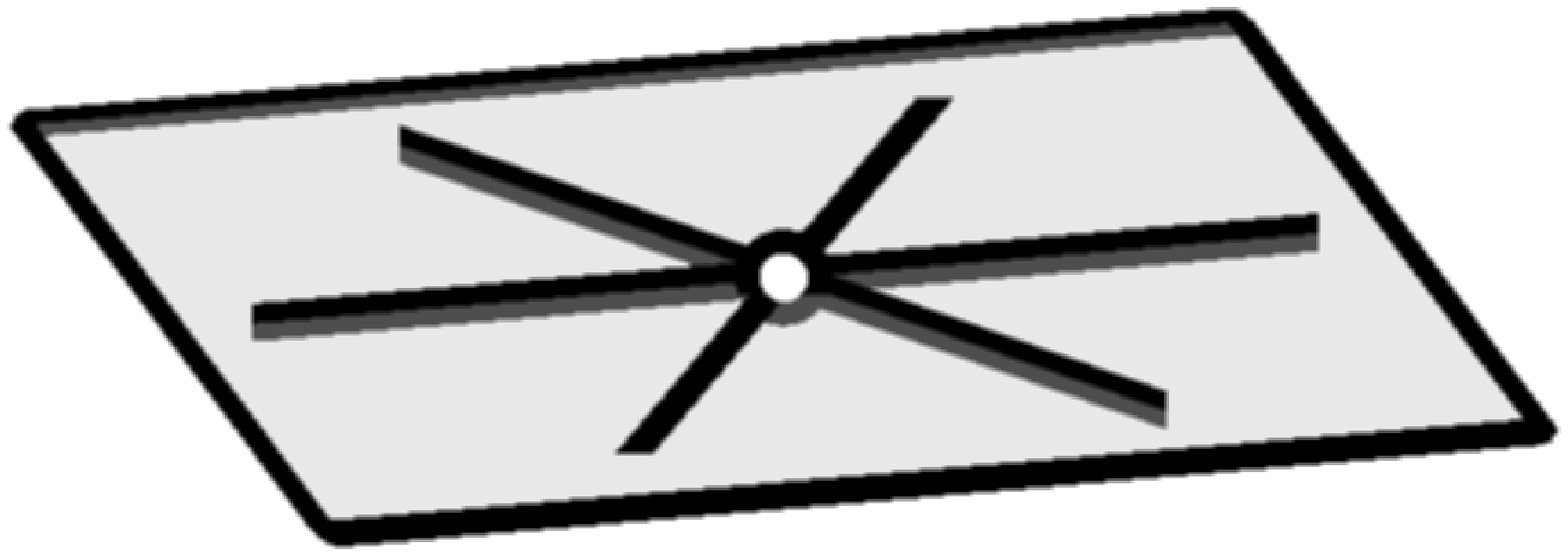}}
         \put(0.85,0.95){$\e$}
         \put(1.95,0.95){$P$}
        \end{picture}
        {\refstepcounter{abbildung}\label{abb3}
          \centerline{Fig.\ \ref{abb3}.}}
      \end{minipage}
      \begin{minipage}[t]{4.0cm}
         \begin{picture}(4.0,4.0)
         \put(0,0)
         {\includegraphics[width=4.0cm]{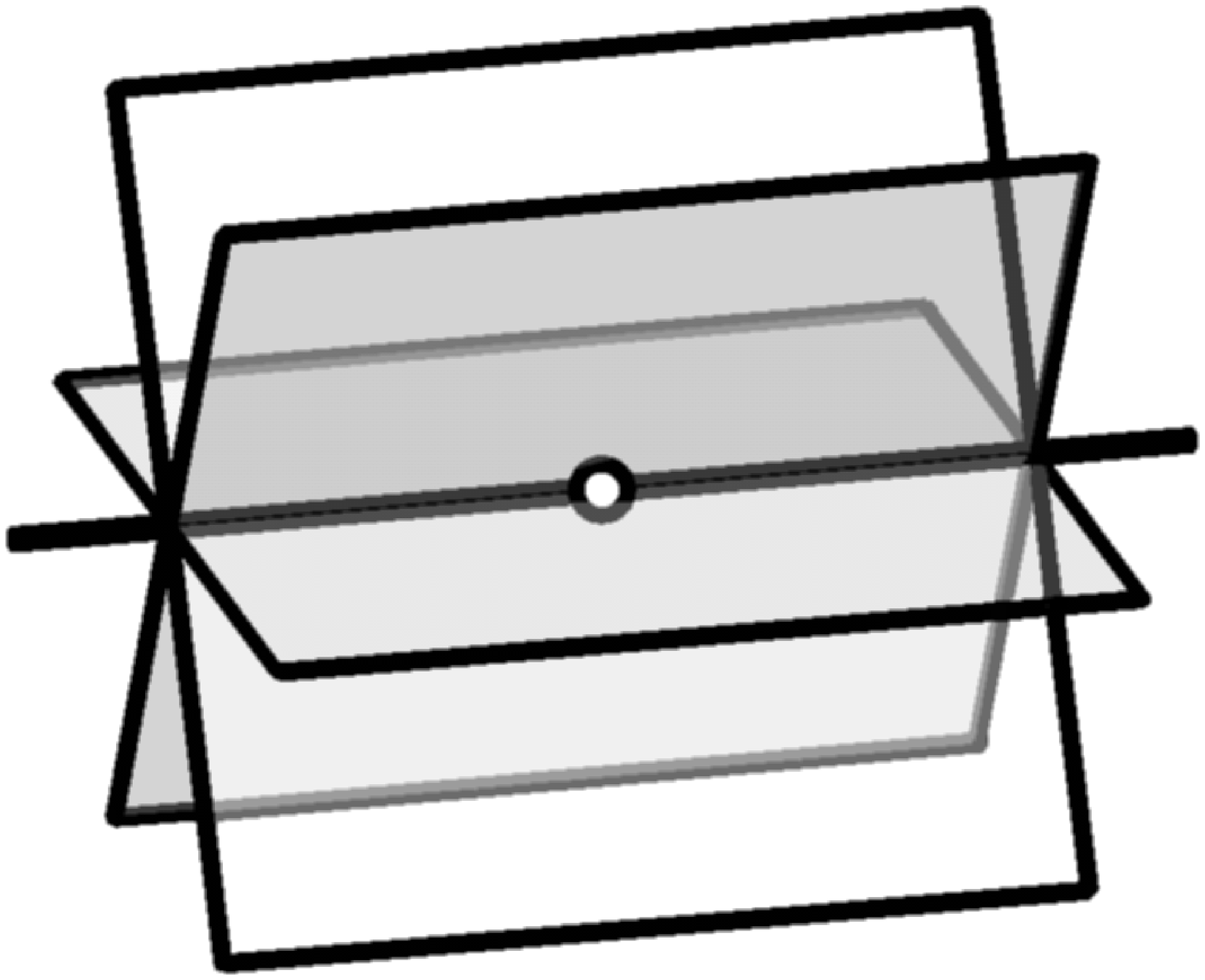}}
         \put(1.6,1.2){$P$}
         \put(3.65,2.0){$g$}
        \end{picture}
        {\refstepcounter{abbildung}\label{abb4}
          \centerline{Fig.\ \ref{abb4}.}}
      \end{minipage}
      \end{center}
}%
If $(P,g,\e)\in \F$ then the pencils given in (\ref{eq:1}) are the
only pencils through it. So each flag is on exactly one pencil of
type $0$, $1$, and $2$. Two distinct pencils of the same type are
disjoint. Two adjacent flags $\Phi, \Psi$ are joined by exactly one
pencil. It will be denoted by $\Phi\Psi$. The following result
describes pencils in the terms of the Pl\"ucker space $(\F, \sim)$:
\begin{pro}\label{pro:1}
  The pencils of flags are exactly the maximal sets of mutually
  related flags.
\end{pro}
\begin{proof}
By definition, the elements of a fixed pencil ${\mathcal M} \subset
\F$ are mutually related. Let $\Phi = (P,g,\e)$ be a flag in
${\mathcal M}$. A flag $\Psi$ is related to $\Phi$ exactly if $\Psi$
is contained in one of the pencils given in (\ref{eq:1}); in
particular one of these pencils is ${\mathcal M}$.

However, if $\Psi \sim \Phi$ is chosen in $\F \bs {\mathcal M}$, then
$\Phi$ is the only element of ${\mathcal M}$ related to $\Psi$,
because every other flag of ${\mathcal M}$ differs from $\Psi$ in
more than one component. So the pencil ${\mathcal M}$ is a maximal
set of mutually related flags.

On the other hand, let ${\mathcal M} \subset \F$ be a maximal set of
mutually related flags. Such an ${\mathcal M}$ contains at least two
adjacent flags, say $\Phi_1$ and $\Phi_2$. We assume that $\Phi_1$
and $\Phi_2$ differ exactly in their $i$-dimensional component. Hence
$\Psi \in {\mathcal M} \bs \{\Phi_1, \Phi_2\}$ implies that the
components of dimension $\ne i$ of $\Psi, \Phi_1$, and $\Phi_2$ are
the same. In other words, $\Psi$ is a flag of the pencil $\Phi_1
\Phi_2$. So ${\mathcal M}\subset \Phi_1 \Phi_2$ and, by the
maximality of ${\mathcal M}$, we have ${\mathcal M} = \Phi_1 \Phi_2$.
\end{proof}
The pair $(\F, \B)$ is a partial linear space \cite[p.~70]{buek-95}
with ``point set'' $\F$ and ``line set'' $\B = \B_0 \cup \B_1 \cup
\B_2$. \textsc{A.~Bichara} and \textsc{C.~Somma} have given an
axiomatic description of this \emph{flag space} associated with the
projective space $(\P,\L)$; see \cite{bich+s-84}, \cite{bich+s-86},
and \cite{bich+s-87}. In terms of this partial linear space related
flags are ``collinear points'' and adjacent flags are ``distinct
collinear points".

\section{Pl\"ucker Transformations}\label{sec:pluecker}
A \emph{Pl\"ucker transformation} of $(\F, \sim)$ is a bijection
$\alpha : \F \to \F$ preserving the relation $\sim$ in both
directions \cite[p.~199]{benz-92}.

From Proposition~\ref{pro:1}, a bijection $\alpha: \F \to \F$ is a
Pl\"ucker transformation exactly if it is an automorphism of the flag
space $(\F , \B)$. Note that here (in contrast to
\cite[p.~61]{bich+s-86}) we do not require that the type of a pencil
is preserved under $\alpha$. It is straightforward to show that each
collineation $\kappa$ of $(\P, \L)$ gives rise to a Pl\"ucker
transformation $\F \to \F: (P,g,\e) \mapsto (P^\kappa, g^\kappa ,
\e^\kappa)$. Similarly, each duality $\delta$ of $(\P, \L)$ yields a
Pl\"ucker transformation $\F \to \F: (P,g,\e) \mapsto (\e^\delta,
g^\delta , P^\delta)$. It is our aim to show that there are no other
Pl\"ucker transformations of $(\F , \sim)$.
\begin{pro}\label{pro:2}
Let $\alpha : \F \to \F$ be a Pl\"ucker transformation of $(\F ,
\sim)$. Then there exists a bijection $\beta : \L \to \L$ such that
  \begin{equation}
    \label{eq:2}
    \F [g]^\alpha = \F [g^\beta] \text{~for all lines~} g \in \L.
  \end{equation}
\end{pro}
\begin{proof}
(a) Choose a line $g \in \L$ and write $\B_i[g]$ for the set of all
pencils of type $i\in \{0,1,2\}$ that are contained in $\F [g]$.
Clearly, $\B_1 [g] = \leer$, whereas $\B_0 [g]$ consists of all
pencils $\F [g,\e]$, $\e \supset g$ an arbitrary plane, and $\B_2
[g]$ consists of all pencils $\F[P,g]$, $P\in g$ an arbitrary point.

Each flag $\Phi \in \F[g]$ is on a unique pencil of $\B_0[g]$ and on
a unique pencil of $\B_2[g]$. Further, each pencil of $\B_0[g]$, say
$\F[g,\e]$, and each pencil of $\B_2[g]$, say $\F [P,g]$, meet at
exactly one flag, namely $(P,g,\e)$. Finally, every pencil has at
least three elements. This means that the incidence structure $(\F
[g], \B_0[g], \B_2[g])$ is a $2$-net \cite[p.~79--80]{kar+k-88}. Cf.\
also \cite[p.~99]{bich+s-84}.

(b) Let, as before, $g\in \L$. We claim that under $\alpha$ no pencil
of $\B_0 [g]$ goes over to a pencil of type 1: Assume to the contrary
that $\F[g, \e]$, $\e \supset g$ a plane, is such a pencil. There
exist distinct points $P,Q \in g$ and a plane $\phi \supset g$ other
than $\e$. We put
 \[ \Phi' := (P,g,\e)^\alpha ,\;
    \Psi' :=(Q,g, \e)^\alpha,\;
    \Phi'':= (P,g,\phi)^\alpha,
    \Psi'':= (Q,g, \phi)^\alpha.
 \]
As $\Phi' \Psi'$ is a pencil of type 1, we get $\Phi'= (P',g',\e')$,
$\Psi' =(P',h', \e')$ with distinct lines $g', h' \in \L$. Since
$\Phi' \Psi'$ is the only pencil of type $1$ through $\Phi'$ and
$\Psi'$, the pencils $\Phi' \Phi''$ and $\Psi' \Psi''$ cannot be of
type 1, whence the line components of $\Phi''$ and $\Psi''$ are $g'$
and $h'$, respectively. However, $g'\cap h'=P'$ and $g' \vee h' =\e'$
implies $\Phi'' \Psi'' = \F[P', \e'] =\Phi' \Psi'$ so that
$\Phi'\sim\Psi''$ which contradicts $(P,g,\e) \not\sim (Q,g, \phi)$.

Similarly, no pencil of $\B_2[g]$ goes over to a pencil of type 1.

(c) Let $(P,g,\e)$ and $(Q,g,\phi)$ be distinct flags. We put
$(P,g,\e)^\alpha =: (P', g' , \e')$. From (a), the pencils $\F[g,\e]
\in \B_0[g]$ and $\F[Q,g] \in \B_2 [g]$ meet at $(Q,g,\e) \in \F[g]$.
Now (b) implies that $\F [g, \e]^\alpha$ and $\F[Q,g]^\alpha$ both
are not of type 1. Hence $g'$ is also line component of
$(Q,g,\e)^\alpha \in \F[g,\e]^\alpha\cap\F[Q,g]^\alpha$ and
$(Q,g,\phi)^\alpha \in \F[Q,g]^\alpha$. Consequently, $\F[g]^\alpha
\subset \F[g']$.

In the same manner we obtain $\F[g']^{\alpha^{-1}} \!\! \subset
\F[g]$, whence $\F[g]^\alpha = \F[g']$. Therefore, by (\ref{eq:2}) we
have a well-defined mapping $\beta: \L \to \L$. Similarly,
$\alpha^{-1}$ defines a mapping $\L\to \L$ which is easily seen to be
the inverse of $\beta$; so $\beta$ is bijective.
\end{proof}
We note that in the terminology of \cite{bich+s-84},
\cite{bich+s-86}, and \cite{bich+s-87} the four pencils $\Phi'\Psi'$,
$\Psi'\Psi''$, $\Psi''\Phi''$, and $\Phi''\Phi'$ that have been
introduced in (b) form a \emph{closed $4$-path}. Such a path cannot
contain a pencil of type $1$. This is one of the axioms used in the
cited papers.

From part (b) above, it is clear that pencils of type~1 go over to
pencils of the same type under $\alpha^{-1}$ as well as under
$\alpha$.

\begin{pro}\label{pro:3}
The bijection $\beta:\L \to \L$ defined in (\ref{eq:2}) and its
inverse mapping $\beta^{-1}$ take intersecting lines to intersecting
lines.
\end{pro}
\begin{proof}
Suppose that $g,h \in \L$ intersect, i.e., $g\cap h =: P$ is a point
and $g \vee h =:\e$ is a plane.

We infer that $\Phi := (P,g,\e) \in \F [g]$ and $\Psi := (P,h,\e) \in
\F[h ]$ span the pencil $\F [P,\e]$  of type~1. Hence its image under
$\alpha$ is again a pencil of type~1. So $\Phi^\alpha = (P', g^\beta
, \e')$ implies $\Psi^\alpha = (P', h^\beta , \e')$. Therefore
$g^\beta$ and  $h^\beta$ intersect.

The proof for $\beta^{-1}$ runs in an analogous way.
\end{proof}
We are now in a position to show the announced result.
\begin{thm}\label{thm:1}
Let $\alpha : \F\to \F$ be a Pl\"ucker transformation of $(\F,
\sim)$. Then there exists either a unique collineation $\kappa$ of
$(\P,\L)$ with
  \begin{equation}
    \label{eq:3}
    (P,g,\e)^\alpha = (P^\kappa, g^\kappa , \e^\kappa)
  \end{equation}
or a unique duality $\delta$ of $(\P, \L)$ with
  \begin{equation}
    \label{eq:4}
    (P,g,\e)^\alpha = (\e^\delta, g^\delta , P^\delta)
  \end{equation}
for all $(P,g,\e)\in\F$.
\end{thm}
  \begin{proof}
From Propositions~\ref{pro:2} and \ref{pro:3}, the given Pl\"ucker
transformation $\alpha$ determines a bijection $\beta:\L \to \L$ such
that $\beta$ and $\beta^{-1}$ map intersecting lines to intersecting
lines. By a result of \textsc{W.L.~Chow} (see
\cite[Theorem~1]{chow-49} or \cite[p.~80--82]{dieu-71}), there exists
either a collineation $\kappa$ of $(\P, \L)$ with $g^\beta =
g^\kappa$ or a duality $\delta$ of $(\P, \L)$ with $g^\beta =
g^\delta$ for all $g\in\L$.

In order to verify (\ref{eq:3}) or (\ref{eq:4}) choose any flag $\Phi
= (P,g,\e) \in \F$. There is a flag $\Psi = (P,h,\e)$ adjacent to
$\Phi$. Now there are two possibilities:

If $\beta$ is induced by a collineation $\kappa$ then $\Phi^\alpha
\in \F [g^\kappa]$  and $\Psi^\alpha \in \F [h^\kappa]$ are adjacent
too. But the only flag in $\F [g^\kappa ]$ that is adjacent to some
flag of $\F[h^\kappa]$ is $(g^\kappa \cap h^\kappa , g^\kappa ,
g^\kappa\vee h^\kappa) = (P^\kappa , g^\kappa, \e^\kappa)$, whence
(\ref{eq:3}) holds.

If $\beta$ is induced by a duality $\delta$ then $\Phi^\alpha =
(g^\delta \cap h^\delta,
g^\delta , g^\delta \vee h^\delta) = (\e^\delta, g^\delta ,
P^\delta)$ follows similarly.

Finally from (\ref{eq:3}) or (\ref{eq:4}), the collineation $\kappa$
or the duality $\delta$ is uniquely determined, since each point
$P\in \P$ is a component of at least one flag.
\end{proof}

\section{The flag variety}\label{sec:variety}

In this section let $(\P,\L)$ be a $3$-dimensional pappian projective
space. We denote by $V$ the underlying vector space with (commutative)
ground field $K$. So $\dim V=4$. Furthermore, the points, lines, and
planes of $(\P,\L)$ are the $1$-, $2$-, and $3$-dimensional subspaces
of $V$, respectively.

Put $(\hatP,\hatL)$ for the projective space on the ($6$-dimensional)
exterior square $V\wedge V$ of $V$. Recall the \emph{Klein mapping}
\begin{equation}\label{eq:5}
  \gamma : \L \to \hatP : \; Kq+Kr \mapsto K(q\wedge r),
\end{equation}
where $q,r\in V$ are linearly independent. It is injective and its
image $\Qua:=\L^\gamma$ is the \emph{Klein quadric} representing the
lines of the projective space $(\P,\L)$. See, for example,
\cite[p.~301--302]{bura-61}, \cite[p.~224]{gier-82}, or
\cite[p.~28--31]{hirs-85}.

Further, let $V^*$ be the dual space of $V$. The $1$-dimensional
subspaces of $V^*$ correspond bijectively to the $3$-dimensional
subspaces of $V$ via $Ke^*\mapsto \ker e^*$ ($e^*\in V^*\bs\{0\}$).
We shall identify the planes of $(\P,\L)$ with the $1$-dimensional
subspaces of $V^*$ or, in other words, the points of the projective
space $(\P^*,\L^*)$ on $V^*$.

Next, we consider the projective space $(\tilP,\tilL)$ on the
($96$-dimensional) tensor product $V\otimes (V\wedge V)\otimes V^*
=:\tilV$. The \emph{Segre mapping}
\begin{equation}\label{eq:6}
  \sigma : \P\times\hatP\times \P^* \to \tilP: \;
  (Kp,Kt,Ke^*) \mapsto K(p\otimes t \otimes e^*),
\end{equation}
where $p\in V$, $t\in V\wedge V$, and $e^*\in V^*$ are non-zero, is
injective and its image is a \emph{Segre variety} $\Se$ of type
$(3,5,3)$ \cite[p.~111]{bura-61}, \cite[Chapter~25.5]{hirs+t-91}. If
we restrict $\sigma$ to the product $\P\times \Qua\times \P^*$ then
we get a point model for all triples consisting of a point, a line,
and a plane of $(\P,\L)$. In particular, from (\ref{eq:5}) and
(\ref{eq:6}) we obtain an injective mapping
\begin{equation}\label{eq:7}
  \phi : \F\to \tilP : (P,g,\e) \mapsto
  (P,g^\gamma,\e)^\sigma
\end{equation}
whose image $\G:=\F^\phi$ is a variety representing the flags of
$(\P, \L)$. The following property of $\G$ is essential.

\begin{pro}\label{pro:4} The $\phi$-images of the pencils of flags are exactly the lines
contained in the flag variety $\G$.
\end{pro}

\begin{proof}
Let $(P,T,\e) \in \P\times\hatP\times\P^*$. Then
  \[(\P \times \{T\} \times \{\e\})^\sigma , (\{P\}\times \hatP \times \{\e\})^\sigma,
  (\{P\}\times \{T\} \times \P^*)^\sigma\]
are the only maximal subspaces contained in the Segre variety $\Se$
that pass through the point $(P,T,\e)^\sigma$
\cite[p.~127--128]{bura-61}. Furthermore, by (\ref{eq:6}), the
mapping $Q \mapsto (Q, T, \e)^\sigma$ is a collineation $\P \to
(\P\times\{T \} \times \{\e\})^\sigma$. Similarly, we have
collineations $\hatP \to (\{P\}\times \hatP \times \{\e\})^\sigma$,
and $\P^* \to (\{P\}\times \{T\} \times \P^*)^\sigma$ (cf.\
\cite[Theorem~25.5.2]{hirs-85}).

Suppose we are given a line $\ell$ contained in $\G$. Choose an
arbitrary point of that line, say $(P,T,\e)^\sigma$. Since $\ell$ is
also a line of the Segre variety $\Se$, its $\sigma$-preimage can be
found with the inverse of one of the three collineations described
above. By taking into account that the lines on the Klein quadric are
exactly the $\gamma$-images of the pencils of lines in $\L$
\cite[p.~301]{bura-61}, we see that $\ell$ is the $\phi$-image of a
pencil of flags.

It is immediately clear from (\ref{eq:5}) and (\ref{eq:6}) that under
$\phi$ each pencil of flags is mapped onto a line contained in $\G$.
\end{proof}

The flag variety $\G$ is the intersection of the Segre variety $\Se$
with a subspace of $(\tilP,\tilL)$ \cite[Satz~2]{bura-54}. More
precisely, we have the following:

\begin{pro}\label{thm:2}
The variety $\G$ representing the flags of
$(\P,\L)$ is the intersection of the Segre variety $\Se$ given by
(\ref{eq:6}) with a $63$-dimensional projective subspace
$(\barP,\barL)$ of the $95$-dimensional projective space
$(\tilP,\tilL)$.
\end{pro}

\begin{proof}
(a) We shall argue in terms of the ($16$-dimensional) exterior
algebra $\Wedge V$, its dual $(\Wedge V)^*$ which will be identified
with $\Wedge V^*$, and the inner product $\iprod:\Wedge V\times
\Wedge V^*\to \Wedge V$. See, among others, \cite{bour-89},
\cite{dehe-93}, or \cite{kowa-79}.

(b) Our first aim is to describe incidence of a point and a line: The
mapping
\begin{equation*}
  V\times (V\wedge V)\times V^* \to \Wedge^3 V\otimes V^* :
  (p,t,e^*)\mapsto (p\wedge t)\otimes e^*
\end{equation*}
is trilinear. By the universal property of the tensor product
$V\otimes (V\wedge V)\otimes V^* = \tilV$, there is
a unique linear mapping
\begin{equation*}
  i_{01} : \tilV \to \Wedge^3 V\otimes V^*
\end{equation*}
with  $p\otimes t\otimes e^*\mapsto (p\wedge t)\otimes e^*$ for all
$(p,t,e^*)\in V\times (V\wedge V)\times V^*$. As $i_{01}$ is
surjective, the dimension of $I_{01}:=\ker i_{01}$ equals $96-16=80$.

Choose any triple $(p,t,e^*)\in V\times (V\wedge V)\times V^*$ with
$p,t,e^*\neq 0$. Then
\begin{equation}\label{eq:8}
  p\otimes t\otimes e^*\in I_{01} \Leftrightarrow  p\wedge t = 0.
\end{equation}
The subspace $T:=\{x\in V \mid x\wedge t=0 \}$ is at most
$2$-dimensional and $\dim T=2$ characterizes $t$ as being
decomposable \cite[47.5]{kowa-79}. Further, the product of the
bilinear mapping $V\times V\to \Wedge^4 V : (v,w)\mapsto v\wedge
w\wedge t$ with an (arbitrarily chosen) isomorphism $\Wedge^4 V \to
K$ is a non-zero alternating bilinear form with radical $T$. The rank
of this form is necessarily even. So, it follows that $\dim
T\in\{0,2\}$. We infer from $p\neq 0$ that the right hand side of
(\ref{eq:8}) is equivalent to the existence
of $q,r\in V$ such that $t=q\wedge r$ and such that the point $Kp$ is
on the line $Kq+Kr$ represented by the point $Kt$ of the Klein
quadric.

(c) Next, we turn to the incidence of a line and a plane: The mapping
\begin{equation*}
  V\times (V\wedge V)\times V^* \to V\otimes V :
  (p,t,e^*)\mapsto p\otimes(t\iprod e^*)
\end{equation*}
is trilinear. Hence, as before, there is a unique linear mapping
\begin{equation*}
  i_{12} : \tilV \to V\otimes V
\end{equation*}
with  $p\otimes t\otimes e^*\mapsto p\otimes(t\iprod e^*)$ for all
$(p,t,e^*)\in V\times (V\wedge V)\times V^*$. The image of $i_{12}$
is the $16$-dimensional tensor product $V\otimes V$, whence
$I_{12}:=\ker i_{12}$ is $80$-dimensional.

Choose any triple $(p,t,e^*)\in V\times (V\wedge V)\times V^*$ with
$p,t,e^*\neq 0$. Then
\begin{equation}\label{eq:9}
  p\otimes t\otimes e^*\in I_{12} \Leftrightarrow  t\iprod e^* = 0.
\end{equation}
The bilinear form $V^*\times V^*\to K : (v^*,w^*)\mapsto \langle t,
v^*\wedge w^*\rangle$ is non-zero and alternating. (Here
$\langle\,,\rangle$ denotes the canonical pairing.)
From \cite[47.4, 47.5]{kowa-79} the rank of this
bilinear form is $2$ exactly if $t$ is decomposable. By the
definition of the inner product,
\begin{equation}\label{eq:10}
   \langle t\iprod e^*,w^*\rangle = \langle t,e^*\wedge w^*\rangle
   \textrm{ for all }w^*\in V^*.
\end{equation}
Suppose that $t\iprod e^*=0$. This implies that $e^*\neq 0$ is in
the radical of the bilinear form
from above so that there are $q,r\in V$ with $t=q\wedge r$. Now
(\ref{eq:10}) gives
\begin{equation}\label{eq:11}
   \langle q\wedge r,e^*\wedge w^*\rangle = \det
   \begin{pmatrix}
     \langle q,e^* \rangle & \langle r,e^* \rangle \\
     \langle q,w^* \rangle & \langle r,w^* \rangle
   \end{pmatrix} = 0
   \textrm{ for all }w^*\in V^*.
\end{equation}
Consequently, $\langle q,e^*\rangle = \langle r,e^* \rangle = 0$. By
reversing these arguments it follows that the right hand side of
(\ref{eq:9}) is equivalent
to the fact that $Kt$ is a point of the Klein quadric which describes
a line of the plane $Ke^*$.

(d) It remains to show that $I_{01}\cap I_{12}$ is a $64$-dimensional
subspace of $\tilV$. We establish instead that $I_{01}+I_{12}=\tilV$
which is equivalent by the dimension formula.

Let $b_0,b_1,b_2,b_3$ be a basis of $V$ and put $b_l^*$ for the
vectors of the dual basis. Then the $96$ product vectors
\begin{equation}\label{eq:12}
  b_i\otimes(b_j\wedge b_k) \otimes b_l^*
\end{equation}
where $i,j,k,l\in\{0,1,2,3\}$ and $j<k$ form a basis of
$\tilV$. This implies that $\tilP$ is
spanned by the $96$ points of $\Se$ that represent all triples formed
by a vertex, an edge, and a face of the tetrahedron
$Kb_0,Kb_1,Kb_2,Kb_3$. Hence also $(\P\times\Qua\times\P^*)^\sigma$
generates $\tilP$. So it is enough to show that the points of the
Segre variety $\Se$ that belong to $I_{01}$ or $I_{12}$ generate a
subspace that contains all points of
$(\P\times\Qua\times\P^*)^\sigma$.

So let $(P,g,\e)\in \P\times\L\times\P^*$ with $P\notin g$ and
$g\not\subset\e$. We distinguish two cases:

If $P\notin\e$ then put $g':=P\vee(g\cap\e)$ and $g'':=(P\vee
g)\cap\e$ (figure \ref{abb5}). As $g$, $g'$, and $g''$ are three
distinct elements of a pencil of lines we obtain that
$(P,g^\gamma,\e)^\sigma$, $(P,{g'}^\gamma,\e)^\sigma\subset I_{01}$,
and $(P,{g''}^\gamma,\e)^\sigma\subset I_{12}$ are three distinct
collinear points of $\Se$, whence $(P,g^\gamma,\e)^\sigma\subset
I_{01}+I_{12}$.

{\unitlength1cm
    \begin{center}
      \begin{picture}(5,3.5)
         \put(0,0)
         {\includegraphics[width=4cm]{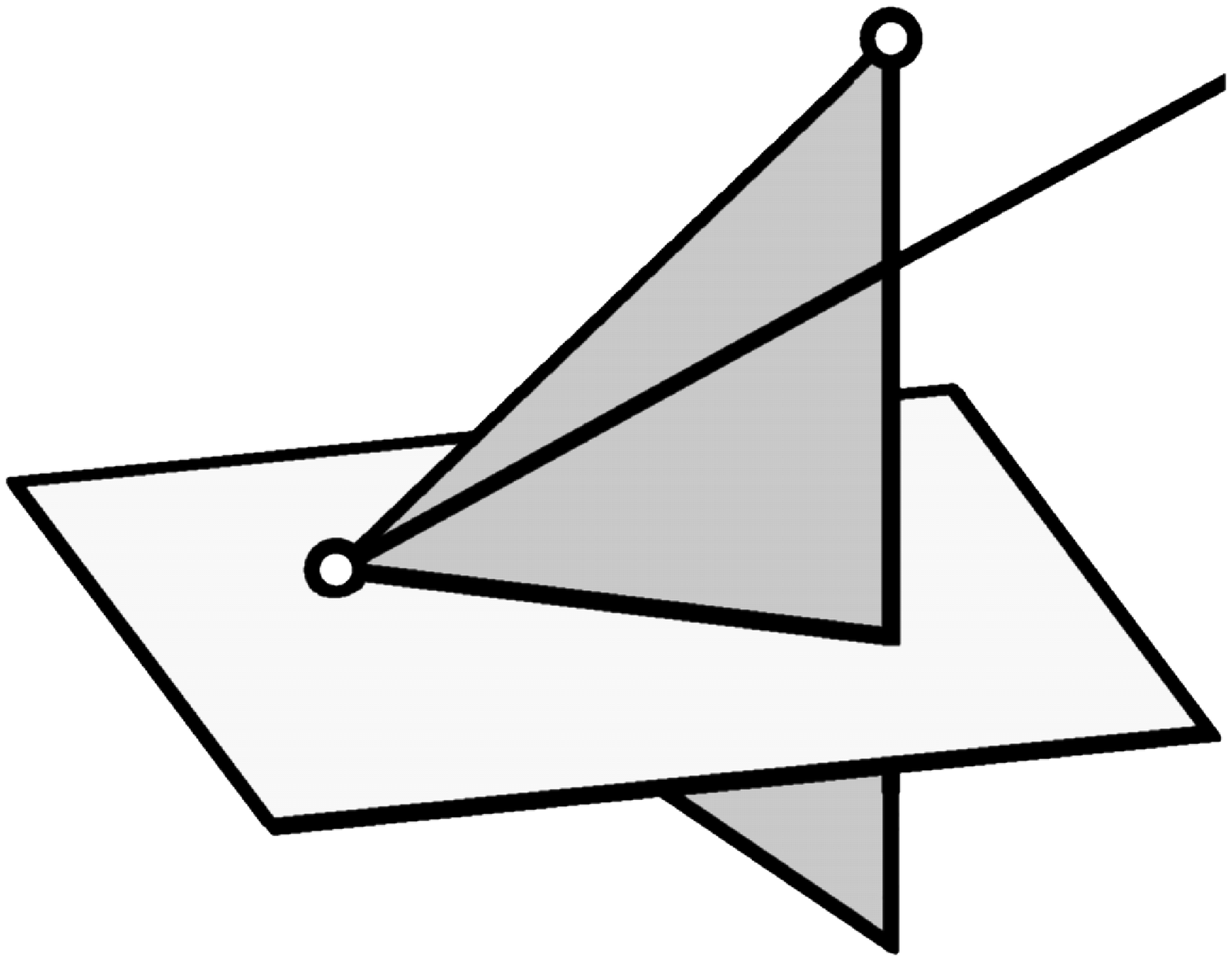}}
         \put(2.8,3.2){$P$}
         \put(3.6,2.4){$g$}
         \put(1.7,2.2){$g'$}
         \put(1.6,.8){$g''$}
         \put(.9,.6){$\e$}
       \end{picture}
       {\refstepcounter{abbildung}\label{abb5}
         \centerline{Fig.\ \ref{abb5}.}}
    \end{center}
 }

If $P\in\e$ then there are points $P',P''\notin \e \cup g$ such that
$P$, $P'$, and $P''$ are three distinct collinear points. We obtain
from the previous case that also $(P,g^\gamma,\e)^\sigma$,
$(P',g^\gamma,\e)^\sigma\subset I_{01}+I_{12}$, and
$(P'',g^\gamma,\e)^\sigma\subset I_{01}+I_{12}$ are three distinct
collinear points of $\Se$, whence $(P,g^\gamma,\e)^\sigma\subset
I_{01}+I_{12}$, as required.
\end{proof}

Clearly, the incidence conditions of (b) and (c) can now easily be
expressed in terms of coordinates; cf.\ also \cite[Satz~1]{bura-54}.

\begin{rem}
The previous result does not answer the question whether or not the
flag variety $\G$ actually spans the $63$-dimensional projective
space $(\barP,\barL)$. In fact, the answer is affirmative. For $K=\C$
this follows from a dimension formula in \cite[p.~142]{bura-58}.
However, at present we can only establish this result for an
arbitrary ground field in terms of $(96\times 64)$-matrices by
explicit computer based calculations (using Maple V). We just sketch
our approach and we use vector space dimensions throughout: Let
$b_0,b_1,b_2,b_3$ be a basis of $V$ and define a basis of $\tilV$ as
in (\ref{eq:12}). Then each flag can be represented by its $96$
homogeneous coordinates with respect to this basis.

In a first step it is easy to show that for each point $Q\in\P$ the
$\phi$-image of $\F[Q]$ spans an $8$-dimensional subspace of $\tilV$;
cf.\ \cite[3.2]{thas+m-99}. Next we consider the four points
$P_i:=Kb_i$ of the coordinate tetrahedron and the four unit points in
the faces of this tetrahedron, i.e. the points
$U_i:=K(b_0+b_1+b_2+b_3-b_i)$ where $i\in\{0,1,2,3\}$. Observe that
the four points $U_i$ are coplanar exactly if the ground field $K$
has characteristic $3$.

The subspace $W_P\subset \tilV$ spanned by
the $\phi$-images of the flags belonging to the union of all subsets
$\F[P_i]$ has dimension $32$. The same holds (irrespective of
$\chara K$) for the subspace $W_U$ spanned by all flags belonging to
the union of all subsets $\F[U_i]$. But now there are two cases:

If $\chara  K\neq 3$ then $\dim (W_P + W_U) = 64$. Otherwise $\dim
(W_P + W_U) = 63$, but the $\phi$-image of the flag given by the
point with coordinates $(1,1,1,-1)$, the line with Pl\"ucker
coordinates $(1,1,-1,0,0,0)$, and the plane with dual coordinates
$(0,0,1,1)$ is not in $W_P + W_U$, whence the assertion follows.
\end{rem}

\begin{rem}
In textbooks on multilinear algebra \emph{Kronecker products} and
\emph{exterior powers} are usually defined only for \emph{linear}
mappings. However, this can easily be extended to semilinear mappings
that share the \emph{same} accompanying automorphism. As an example
we treat the exterior square of a semilinear mapping:

Let $f:X\to Y$ be a semilinear mapping of vector spaces over $K$ with
accompanying automorphism $\zeta\in\Aut(K)$. We define $Y_\zeta$ as
the vector space with the same additive group as $Y$, but with the
modified multiplication $k*y:=k^\zeta y$ for all $k\in K$ and all
$y\in Y$; cf.\ \cite[p.~221]{bour-89}. Then the linear mappings of
$X$ into $Y_\zeta$ are exactly the $\zeta$-semilinear mappings $X\to
Y$. Moreover, $(Y\wedge Y)_\zeta = Y_\zeta\wedge Y_\zeta$. The usual
exterior square of the linear mapping $f:X\to Y_\zeta$ is a linear
mapping $\hat f: X\wedge X\to (Y\wedge Y)_\zeta$ and at the same time
a $\zeta$-semilinear mapping $X\wedge X \to Y\wedge Y$.
\end{rem}

Let us say that two points of the flag variety $\G$ are \emph{related}
if they are on a line which is contained in $\G$. Here is our final
result.

\begin{thm}\label{thm:3}
Let $\eta:\G\to\G$ be a bijection of the variety representing the
flags of a $3$-dimensional pappian projective space such that under
$\eta$ and $\eta^{-1}$ related points go over to related points. Then
$\eta$ extends to a unique collineation of the subspace spanned by
$\G$.
\end{thm}
\begin{proof}
(a) From Proposition \ref{pro:4}, the given bijection $\eta$ is the
$\phi$-transform of a Pl\"ucker transformation $\alpha$ of
$(\F,\sim)$. By Theorem \ref{thm:1}, we obtain that there is either
a collineation $\kappa$ or a duality $\delta$ whose action on $\F$
coincides with $\alpha$.

(b) Let $\kappa$ be such a collineation. Then $\kappa$ is induced by
a semilinear bijection $f:V\to V$. The exterior square of $f$, say
$\hat{f}$, describes the action (via the Klein mapping $\gamma$) of
$\kappa$ on the Klein quadric, and the inverse of the transpose of
$f$, say $f^*: V^* \to V^*$, describes the action of $\kappa$ on the
set of planes. Observe that all three mappings belong to the same
automorphism of $K$. Then their Kronecker product $\tilde f
:= f\otimes \hat f \otimes f^*$ is a semilinear bijection too and
hence induces a collineation $\mu$ of the projective space
$(\tilP,\tilL)$.

(c) Let $\delta$ be such a duality. Then $\delta$ is induced by a
semilinear bijection $f:V\to V^*$. The polarity of the Klein quadric
determines a linear bijection $d:V\wedge V\to V^*\wedge V^*$. (Here
we identify $V^*\wedge V^*$ with the dual of $V\wedge V$). The
product of the exterior square of $f$, say $\hat{f}$,  with $d^{-1}$
describes the action (via the Klein mapping $\gamma$) of $\delta$ on
the Klein quadric, and the inverse of the transpose of $f$, say
$f^*:V^*\to V$, describes the action of $\delta$ on the set of
planes. Now the remaining proof runs as before by virtue of the
commutativity of the tensor product, i.e.\ the canonical isomorphism
$V^*\otimes (V\wedge V)\otimes V \cong V\otimes (V\wedge V)\otimes
V^*$.

(d) The collineation $\mu$ leaves invariant the subspace generated by
the flag variety $\G$ and, by construction, extends the given
bijection $\eta$. On the other hand, let $\rho$ be a
collineation with the required properties. Then $\mu$ and $\rho$
coincide for all lines contained in $\G$. Since $(\F,\sim)$ is
connected, any two points of $\G$ can be joined by a polygonal path
contained in $\G$ consisting of $m$ lines, say. Then it is an easy
induction on $m\geq 1$ that $\mu$ and $\rho$ coincide on the subspace
spanned by the lines of the polygon. Thus, finally, the two
collineations are the same on the subspace spanned by $\G$.
\end{proof}


Addresses of the Authors:\\
~\\
\noindent
 Hans Havlicek, Klaus List, Institut f\"ur Geometrie,
 Technische Universit\"at, Wiedner Hauptstra{\ss}e 8--10, A 1040 Wien,
 Austria\\
 havlicek@geometrie.tuwien.ac.at\\
 klaus@geometrie.tuwien.ac.at

\noindent
 Corrado Zanella, Dipartimento di Matematica Pura ed
 Applicata, Universit\`a di Padova, via Belzoni 7, I 35131 Padova,
 Italy\\
 zanella@math.unipd.it
\end{document}